\documentclass[11pt]{amsart}
\usepackage{amsmath}
\usepackage{amsfonts}
\usepackage{latexsym}
\usepackage{amssymb}
\usepackage{enumerate}
\usepackage[dvips]{graphics}

\newcommand{\Z}{{\mathbf Z}}

\newcommand{\R}{{\mathbf R}}

\newcommand{\comment}[1]{}

\def\A{\mathcal A}

\def\R{\mathbf R}

\newtheorem{theorem}{Theorem}%[section]

\newtheorem{prop}{Proposition}
\newtheorem{lemma}{Lemma}

\newtheorem{definition}{Definition}

\numberwithin{equation}{section}
\begin{document}
\title[Smooth Lyapunov 1-forms]{Smooth Lyapunov 1-forms}
\author[M. ~Farber]{M.~Farber}
\address{Department of Mathematics, Tel Aviv University, Tel Aviv 69978, Israel}
\email{mfarber@tau.ac.il}
\thanks{M. Farber was partially supported by a grant from the Israel Academy of Sciences
and Humanities; this work was done while M. Farber visited FIM ETH in Zurich}

\author[T.~Kappeler]{T. ~Kappeler}
\address{Institute of Mathematics, University of Z\"urich,
8057 Z\"urich, Switzerland}

\email{tk@math.unizh.ch}
\thanks{T. Kappeler and J. Latschev were partially
supported by the European Commission under grant HPRN-CT-1999-00118}

\author[J.~Latschev]{J.~Latschev}
\address{Institute of Mathematics, University of Z\"urich,
8057 Z\"urich, Switzerland} \email{janko@math.unizh.ch}

\author[E.~Zehnder]{E.~Zehnder}
\address{Department of Mathematics, ETH Z\"urich, 8092 Z\"urich, Switzerland}
\email{eduard.zehnder@math.ethz.ch}


\date{February 27, 2003}


\begin{abstract}
We find conditions which guarantee that a given flow $\Phi$ on a closed
smooth manifold $M$ admits a smooth Lyapunov one-form $\omega$ lying in a prescribed de Rham
cohomology class $\xi\in H^1(M;\R)$. These conditions are formulated in terms of
Schwartzman's asymptotic cycles $\A_\mu(\Phi)\in H_1(M;\R)$ of the flow.
\end{abstract}
\maketitle

\sloppy

\section{Introduction}

C. Conley \cite{Conl, Conl1} showed  that any continuous flow $\Phi: X\times
\R\to X$ on a compact metric space $X$ {\it \lq\lq decomposes\rq\rq} into a
chain recurrent flow and a gradient-like flow. More precisely, he proved the
existence of a continuous function $L: X\to \R$ which (i) decreases along any
orbit of the flow in the complement $X-R$ of the chain recurrent set $R\subset
X$ of $\Phi$ and (ii) is constant on the connected components of $R$. Such a
function $L$ is called a {\it Lyapunov function} for $\Phi$. This existence
result plays a fundamental role in Conley's program of understanding general
flows as collections of isolated invariant sets linked by heteroclinic orbits.

A more general notion of a {\it Lyapunov 1-form} was introduced in paper \cite{Fa:2}.
Lyapunov 1-forms, compared to Lyapunov functions,
allow to go one step further and to analyze the
flow within the chain-recurrent set $R$ as well.
Lyapunov 1-forms provide an important tool in applying methods of homotopy theory to dynamical
systems. In the recent papers
\cite{Fa:1}, \cite{Fa:2} a generalization of the Lusternik - Schnirelman theory was constructed which applies to flows
admitting Lyapunov 1-forms.

The problem of existence of Lyapunov 1-forms was addressed in our recent preprint \cite{FKLZ}, where we
worked in the category of compact metric spaces, continuous flows and continuous closed 1-forms.
In the present paper we
study the smooth version of the problem: we construct smooth
Lyapunov 1-forms for smooth flows on smooth manifolds.
We use Schwartzman's asymptotic cycles to formulate a necessary condition for the existence of Lyapunov 1-forms in a given cohomology class.
We also show that under an additional assumption this condition is equivalent to the
homological condition introduced in our previous paper \cite{FKLZ}.

 \section{Definition}

Let $V$ be a smooth vector field on a smooth manifold $M$.
Assume that $V$ generates a continuous flow $\Phi: M\times \R\to M$ and $Y\subset M$ is a closed,
flow-invariant subset.

\begin{definition}\label{def1}
\label{introlyap}
A smooth closed 1-form $\omega$ on $M$ is called a Lyapunov one-form for the pair
$(\Phi, Y)$ if it has the following properties:
\begin{enumerate}
\item[{\rm ($\Lambda$1)}] The function $\iota_V(\omega)=\omega(V)$ is negative on $M-Y$;

\item[{\rm ($\Lambda$2)}] There exists a smooth function $f:U\to \R$ defined on an open neighborhood $U$ of $Y$
such that
$$\omega|_U=df \quad \text{and}\quad df|_Y=0.$$
\end{enumerate}
\end{definition}

The above definition is a modification of the notion of a Lyapunov 1-form introduced in
section 6 of \cite{Fa:2}. The definition of \cite{Fa:2} requires that $Y$ consists of
finitely many points and the vector field $V$ is locally a gradient of $\omega$ with
respect to a Riemannian metric.

Definition  \ref{def1} can also be compared with the definition of a Lyapunov 1-form in
the continuous setting which was introduced in \cite{FKLZ}. Condition ($\Lambda$1) above
is slightly stronger than condition (L1) of Definition 1 in \cite{FKLZ}. Condition
($\Lambda$2) is similar to condition (L2) of Definition 1 from \cite{FKLZ} although they
are not equivalent.

There are several natural alternatives for condition ($\Lambda$2). One of them is:
\noindent
\begin{enumerate}
\item[{\rm ($\Lambda 2^\prime$)}] The 1-form $\omega$, viewed as a map $\omega: M\to T^\ast(M)$, vanishes on $Y$.
\end{enumerate}
\noindent
It is clear that ($\Lambda 2$) implies ($\Lambda 2^\prime$).
We can show that the converse is true under some additional assumptions:

\begin{lemma} If the de Rham cohomology class $\xi$ of $\omega$ is integral, $\xi=[\omega]\in H^1(M;\Z)$, then the conditions
($\Lambda 2^\prime$) and ($\Lambda 2$) are equivalent.
\end{lemma}
\begin{proof}
Clearly we only need to show that ($\Lambda 2^\prime$) implies ($\Lambda 2$).
Since $\xi$ is integral there exists a smooth map $\phi: M\to S^1$ such that
$\omega=\phi^\ast(d\theta)$, where $d\theta$ is the standard angular 1-form on the circle $S^1$.
Let $\alpha\in S^1$ be a regular value of $\phi$. Assuming that ($\Lambda 2^\prime$) holds it then follows that
$U=M - \phi^{-1}(\alpha)$ is an open neighborhood of $Y$. Clearly $\omega|_U=df$ where $f:U \to \R$ is a smooth function which
is related to $\phi$ by $\phi(x)=\exp(if(x))$ for any $x\in U$. Hence ($\Lambda 2$) holds. 
\end{proof}

\begin{lemma} The conditions
($\Lambda 2^\prime$) and ($\Lambda 2$) are equivalent if $Y$ is an Euclidean Neighborhood Retract (ENR).
\end{lemma}
\begin{proof} Again, we only have to establish ($\Lambda 2^\prime$) $\Rightarrow$
($\Lambda 2$). Since $Y$ is an ENR it admits an open neighbourhood $U\subset M$ such that the inclusion $i_U: U\to M$ is homotopic to
$i_Y\circ r$, where $i_Y:Y\to M$ is the inclusion and $r: U\to Y$ is a retraction (see \cite{D},
chapter 4, \S 8, Corollary 8.7). Pick a base point $x_j$ in every path-connected component $U_j$ of $U$ and
define
a smooth function $f_j:U_j\to \R$ by
$$f_j(x) =\int\limits_{x_j}^x \omega, \quad x\in U_j.$$
The latter integral is independent of the choice of the integration path in $U_j$ connecting $x_j$ with $x$.
This claim is equivalent to the vanishing of the integral $\int_\gamma \omega$ for any closed loop $\gamma$ lying in $U$.
To show this we apply the retraction to see that
$\gamma$ is homotopic in $M$ to the loop $\gamma_1=r\circ \gamma$, which lies in $Y$; thus we obtain
$\int_\gamma \omega = \int_{\gamma_1} \omega=0$ because of ($\Lambda 2^\prime$).
It is clear that the functions $f_j$ together determine a smooth function $f:U\to \R$ with $df=\omega|_U$.
\end{proof}

Our main goal in this paper is to find topological conditions which guarantee
that for a given vector field $V$ on $M$ there
exists a Lyapunov 1-form $\omega$ lying in a prescribed cohomology class
$\xi\in H^1(M;\R)$.

\section{Asymptotic cycles of Schwartzman}

Let $M$ be a closed smooth manifold and let $V$ be a smooth vector field.
Let $\Phi: M\times \R\to M$ be the flow generated by $V$.

Consider a Borel measure $\mu$ on $M$ which is invariant under $\Phi$.
According to S. Schwartzman \cite{Sc}, these data determine a real homology class
$$\mathcal A_\mu= \A_\mu(\Phi)\, \in H_1(M;\R)$$
called {\it the asymptotic cycle of the flow $\Phi$ corresponding to the measure $\mu$.}
The class $\A_\mu$ is defined as follows.
For a de Rham cohomology class $\xi \in H^1(M;\R)$ the evaluation $\langle \xi, \A_\mu\rangle \in \R$ is
given by the integral
\begin{eqnarray}
\langle \xi , \A_\mu\rangle  \, =\, \int_M \iota_V(\omega)d\mu,\label{integral}
\end{eqnarray}
where $\omega$ is a closed 1-form in the class $\xi$.
Note that
$\langle \xi, \A_\mu\rangle$ is well-defined, i.e. it depends only on the cohomology class $\xi$ of $\omega$,
see \cite{Sc}, page 277.
Indeed, replacing $\omega$ by $\omega'=\omega+df$, where $f: M\to \R$ is a smooth function,
the integral in (\ref{integral}) gets changed by the quantity
\begin{eqnarray}
\int_M V(f)d\mu = \lim_{s\rightarrow 0} \;\frac 1 s \int_M \bigl\{ f(x\cdot s) -
f(x)\bigr\}\, d\mu(x) .\label{integral1}
\end{eqnarray}
Here $V(f)$ denotes the derivative of $f$ in the direction of the vector field $V$ and
$x\cdot s$ stands for the flow $\Phi(x, s)$ of the vector field $V$. Since the measure $\mu$ is flow invariant, the integral
on the RHS of (\ref{integral1}) vanishes for any $f$. It is clear that the RHS of
(\ref{integral}) is a linear function of $\xi\in H^1(M;\R)$. Hence there exists a unique
real homology class $\A_\mu\in H_1(M;\R)$ which satisfies (\ref{integral}) for all
$\xi\in H^1(M;\R)$.

\section{Necessary Conditions}

We consider the flow $\Phi$ as being fixed and we vary the invariant measure $\mu$.
As the class $\A_\mu\in H_1(M;\R)$ depends linearly on $\mu$, the set of asymptotic cycles $\A_\mu$ corresponding to all
$\Phi$-invariant positive measures $\mu$ forms a convex cone in the vector space $H_1(M;\R)$.

\begin{prop} Assume that there exists a Lyapunov 1-form for $(\Phi,Y)$ lying in a cohomology class
$\xi\in H^1(M;\R)$. Then
\begin{eqnarray}
\langle \xi, \A_\mu\rangle \, \leq \, 0\label{new11}
\end{eqnarray}
for any $\Phi$-invariant positive Borel measure $\mu$ on $M$; equality in
(\ref{new11})
takes place if and only if
the complement of $Y$ has measure zero.
Further, the restriction of $\xi$ to $Y$, viewed as a \v Cech cohomology class
$$\xi|_{Y}\in \check H^1(Y;\R)$$
vanishes, $\xi|_{Y} =0$.
\end{prop}
\begin{proof} Let $\omega$ be a Lyapunov one-form for $(\Phi, Y)$ lying in the class $\xi$.
According to Definition \ref{def1}, the function $\iota_V(\omega)$ is negative on $M-Y$ and vanishes on $Y$. We obtain that the integral
$$\int_M\iota_V(\omega)d\mu = \langle \xi, \A_\mu\rangle$$
is nonpositive.

Assuming $\mu(M-Y)>0$, we find a compact $K\subset M-Y$ with $\mu(K)>0$; this follows from
the Theorem of Riesz - see e.g. \cite{La}, Theorem 2.3(iv), page 256.
There is a constant $\epsilon >0$ such that $\iota_V(\omega)|_K \leq -\epsilon$. Therefore,
one has
$$\int_M \iota_V(\omega)d\mu \leq -\epsilon \mu(K) \, <\, 0.$$
Hence, the value $\langle\xi,\A_\mu\rangle$ is strictly negative if the measure $\mu$ is not supported in $Y$.

To prove the second statement we observe (see \cite{Sp}) that the \v Cech cohomology $\check H^1(Y;\R)$ equals the direct limit of the singular cohomology
$$\check H^1(Y;\R) = \lim_{W\supset Y} H^1(W;\R),$$
where $W$ runs over open neighborhoods of $Y$. It is clear in view of condition ($\Lambda$2) that $\xi|_U=0\in H^1(U;\R)$
(by the de Rham theorem). Hence the result follows.
\end{proof}

\section{Chain-recurrent set $R_\xi$}

Given a flow $\Phi$, our aim is to construct a Lyapunov 1-form $\omega$ for a pair $(\Phi, Y)$
lying in a given cohomology class
$\xi\in H^1(M;\R)$.
A natural candidate for $Y$ is the subset $R_\xi=R_\xi(\Phi)$ of the chain-recurrent set $R=R(\Phi)$ which was defined in
\cite{FKLZ}. For convenience of the reader we briefly recall the definition.

Fix a Riemannian metric on $M$ and denote by $d$ the corresponding distance function.
Given any $\delta>0$, $T>1$, a $(\delta,T)$-{\it chain from $x\in M$ to $y\in
M$}  is a finite sequence $x_0=x, x_1, \dots, x_{N}=y$ of points in $M$
and numbers $t_1, \dots, t_N\in \R$ such that $t_i\geq T$ and
$d(x_{i-1}\cdot t_i, x_i) < \delta$ for all $1 \leq i \leq N$. Here we use the notation $\Phi(x,t)=x\cdot t$.
The {\it chain recurrent set} $R=R(\Phi)$ of the flow
$\Phi$ is defined as the set of all points $x\in M$ such that for any
$\delta>0$ and $T > 1$ there exists a $(\delta,T)$-chain starting and ending
at $x$. The chain recurrent set is
closed and invariant under the flow.

Given a cohomology class $\xi\in H^1(M;\R)$ there is a natural covering space $p_\xi: \tilde M_\xi\to M$ associated with $\xi$.
A closed loop $\gamma: [0,1]\to M$ lifts to a closed loop in $\tilde M_\xi$
if and only if the value of the cohomology class $\xi$ on the homology class $[\gamma]\in H_1(M;\Z)$ vanishes,
$\langle \xi, [\gamma]\rangle =0$. See \cite{Sp}.

The flow $\Phi$ lifts uniquely to a flow $\tilde \Phi$ on the covering $\tilde M_\xi$.
Consider the chain recurrent set $R(\tilde \Phi)\subset \tilde M_\xi$
of the lifted flow and denote by $R_\xi=p_\xi(R(\tilde \Phi))\subset M$ its projection onto
$M$.
The set $R_\xi$ is referred to as {\it the chain recurrent set associated to the cohomology class $\xi$.}
It is clear that $R_\xi$ is a closed and $\Phi$-invariant subset of $R$. We denote by $C_\xi$ the complement of $R_\xi$ in $R$,
$$C_\xi = R- R_\xi.$$

A different definition of $R_\xi$ which does not use the covering space $\tilde M_\xi$ can be found in
\cite{FKLZ}.

To state our main result we also need the following notion.
\begin{figure}[h]
\begin{center}
%\resizebox{9cm}{5cm}{\includegraphics[155,540][440,675]{knot1.eps}}
\resizebox{7cm}{4cm} {\includegraphics[105,448][437,659]{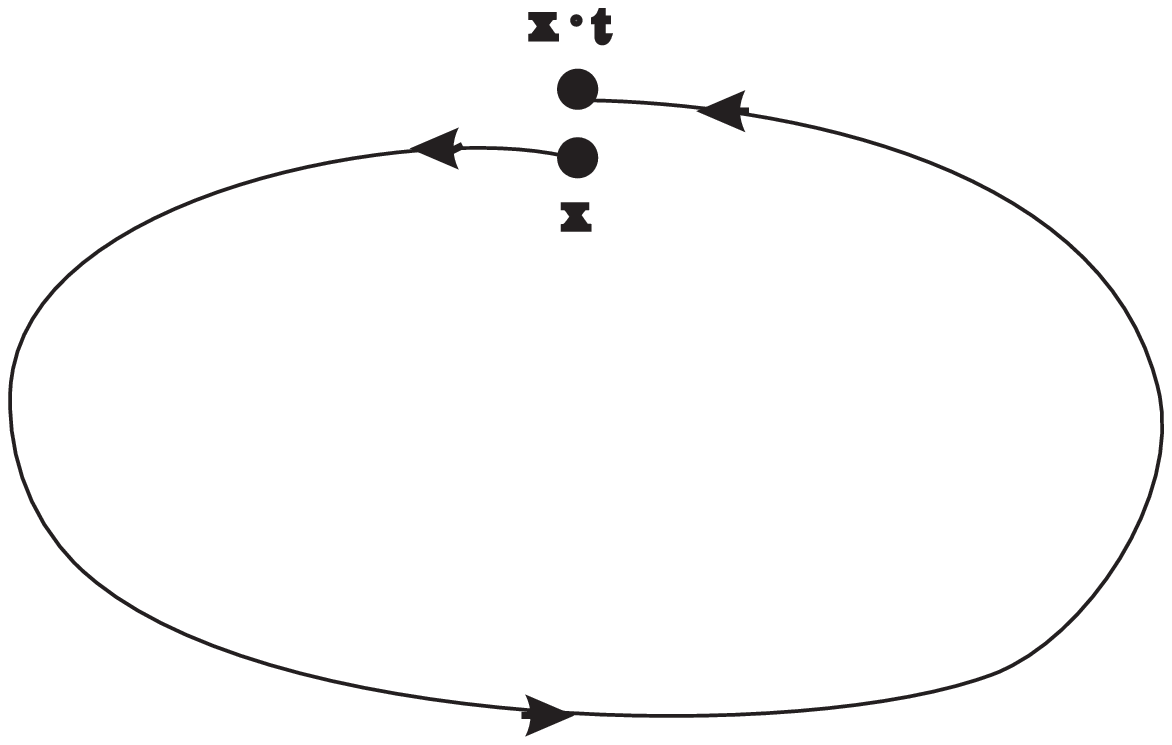}}
\end{center}
\end{figure}
A $(\delta, T)${\it -cycle} of the flow $\Phi$ is defined as a pair $(x,t)$, where $x\in
M$ and $t>T$ such that $d(x, x\cdot t)<\delta$. If $\delta$ is small enough then any
$(\delta, T)$-cycle determines in a canonical way a unique homology class $z\in
H_1(M;\Z)$ which is represented by the flow trajectory from $x$ to $x\cdot t$ followed by
a \lq\lq {\it short }\rq\rq arc connecting $x\cdot t$ with $x$. See \cite{FKLZ}.

\section{Theorem}

\begin{theorem}\label{main}
Let $V$ be a smooth vector field on a smooth closed manifold $M$. Denote by
$\Phi: M\times \R\to M$ the flow generated by $V$.
Let $\xi\in H^1(M;\R)$ be a cohomology class such that the restriction $\xi|_{R_\xi}$, viewed as a \v Cech cohomology class
$\xi|_{R_\xi}\in \check H^1(R_\xi;\R)$,
vanishes. Then the following properties of $\xi$ are equivalent:

{\rm (I)}. There exists a smooth Lyapunov 1-form for $(\Phi, R_\xi)$ in the cohomology
class $\xi$ and the subset $C_\xi$ is closed.

{\rm (II)}. For any Riemannian metric on $M$ there exist $\delta>0$ and $T>1$ such that
the homology class $z\in H_1(M;\Z)$ associated with an arbitrary $(\delta, T)$-cycle
$(x,t)$ of the flow, with $x\in C_\xi$, satisfies $\langle \xi, z\rangle \leq -1$.

{\rm (III)}. The subset $C_\xi$ is closed and there exists a constant $\eta>0$
such that for any $\Phi$-invariant positive Borel measure $\mu$ on $M$ the asymptotic
cycle $\A_\mu=\A_\mu(\Phi)\in H_1(M;\R)$ satisfies
\begin{eqnarray}
\langle \xi, \A_\mu\rangle \, \leq\,   -\eta\cdot \mu(C_\xi).\label{inequa3}
\end{eqnarray}

{\rm (IV)}. The subset $C_\xi$ is closed and for any $\Phi$-invariant positive Borel
measure $\mu$ on $X$ with $\mu(C_\xi)>0$, the asymptotic cycle $\A_\mu=\A_\mu(\Phi)\in
H_1(M;\R)$ satisfies
\begin{eqnarray}
\langle \xi, \A_\mu\rangle < 0.\label{inequa}
\end{eqnarray}
\end{theorem}

The main point of this result is that it
gives sufficient homological conditions for the existence of a Lyapunov 1-form in the cohomology class $\xi$.

Condition (\ref{inequa3}) can be reformulated using the notion of a quasi-regular point.
Recall that $x\in X$ is a {\it quasi-regular} point of the flow $\Phi: X\times \R\to X$
if for any continuous function $f: X\to \R$ the limit
\begin{eqnarray}
\lim_{t\to \infty} \frac{1}{t} \int_0^t f(x\cdot s)ds
\end{eqnarray}
exists. It follows from the ergodic theorem that the subset $Q\subset X$ of all
quasi-regular points has full measure with respect to any $\Phi$-invariant positive Borel
measure on $X$, see \cite{Ja}, page 106. From the Riesz representation theorem, see e.g.
\cite{Ru}, page 256, one deduces that for any quasi-regular point $x\in Q$ there exists a
unique positive flow-invariant Borel measure $\mu_x$ with $\mu_x(X)=1$ satisfying
\begin{eqnarray}
\lim_{t\to \infty} \frac{1}{t} \int_0^t f(x\cdot s)ds = \int_X fd\mu_x
\end{eqnarray}
for any continuous function $f$. We use below the
well-known fact that any positive,
$\Phi$-invariant Borel measure $\mu$ with $\mu(X)=1$ belongs to the weak$^\ast$ closure
of the convex hull of the set of measures $\mu_x$, $x\in Q$, see \cite{Ja}, p. 108.

If the subset $C_\xi\subset X$ is closed, and hence compact, one can apply the above
mentioned facts to the restriction of the flow to $C_\xi$.
Let $\omega$ be an arbitrary
smooth closed 1-form lying in the cohomology class $\xi$. For any quasi-regular point
$x\in C_\xi$ of the flow $\Phi|_{C_\xi}$ one has
\begin{eqnarray}
\lim_{t\to \infty} \frac{1}{t}\int_x^{x\cdot t} \omega = \lim_{t\to \infty}\frac{1}{t}
\int_0^t \iota_V(\omega)(x\cdot s)ds=\\
= \int_{M}\iota_V(\omega)d\mu_x = \langle \xi, \mathcal A_{\mu_x}\rangle.\nonumber
\end{eqnarray}

We therefore conclude that condition {\rm (III)} is equivalent to:

{$\rm (III')$} {\it The subset $C_\xi$ is closed and there exists a constant $\eta>0$ such that for any quasi-regular point $x\in
C_\xi$,
\begin{eqnarray}
\lim_{t\to \infty} \frac{1}{t}\int_x^{x\cdot t} \omega\, \, \leq\, \, -\eta,\label{stam2}
\end{eqnarray}
where $\omega$ is an arbitrary closed 1-form in class $\xi$.}

The value of the limit (\ref{stam2}) is independent of the choice of a closed 1-form
$\omega$; the only requirement is that $\omega$ lies in the cohomology class $\xi$.

In the special case $\xi=0$ the set $C_\xi$ is empty and $R=R_\xi$. The above
statement then reduces to the following well-known theorem of C.Conley - see
\cite{Conl} and  \cite{Sh}, Theorem 3.14:

\begin{prop}\label{conley} {\rm (C. Conley)} Let $V$ be a smooth vector field on a smooth closed manifold $M$. Denote by
$\Phi: M\times \R\to M$ the flow generated by $V$ and by $R$ the chain recurrent set of
$\Phi$. Then there exists a smooth Lyapunov function $L: M\to \R$ for $(\Phi, R)$. This
means that $V(L)<0$ on $M-R$ and $dL=0$ pointwise on $R$.
\end{prop}

Proposition \ref{conley} is used in the proof of Theorem \ref{main}.

As we could not find a proof of this statement in the literature we present one in the
appendix.

\section{Proof of Theorem \ref{main}}

The implication (I) $\Rightarrow $ (II) follows from the proof of Proposition 4 in \cite{FKLZ}.

(II) $\Rightarrow$ (III) By \cite{FKLZ}, Theorem 2, the set $C_\xi$ is closed. Now we want to show that the inequality
(\ref{inequa}) is satisfied for any positive $\Phi$-invariant Borel measure $\mu$ on $X$ with $\mu(C_\xi)>0$.
Fix a closed 1-form $\omega$ in the cohomology class $\xi$. By Lemma 6 from \cite{FKLZ}, there exist constants $\alpha>0$ and $\beta>0$
such that for any $x\in C_\xi$ and $t>0$, one has
$$\int_x^{x\cdot t}\omega \leq -\alpha t +\beta.$$
Set $t_0={2\beta}/\alpha$. Then for any $x\in C_\xi$ and $t\geq t_0$ we have
\begin{eqnarray}
\frac{1}{t} \int_x^{x\cdot t} \omega \leq -\frac{\alpha}{2}.\label{bound}
\end{eqnarray}
With any quasi-regular point $x\in C_\xi$ one associates in a canonical way a
positive $\Phi$-invariant Borel
measure $\mu_x$ on $C_\xi$, see above.
It has the property that
\begin{eqnarray}
\lim_{t\to \infty} \frac{1}{t}\int_x^{x\cdot t}\omega = \int_M \iota_V(\omega)d\mu_x.\label{stam}
\end{eqnarray}
From (\ref{bound}) and (\ref{stam}) one obtains
\begin{eqnarray}
\langle \xi, \A_{\mu_x}\rangle \leq -\frac{\alpha}{2} <0
\end{eqnarray}
for any quasi-regular point $x\in C_\xi$.
According to \cite{Ja}, page 108,
any positive $\Phi$-invariant Borel measure $\mu$ with $\mu(M) = \mu(C_\xi)=1$ belongs to the weak$^\ast$ closure
of the convex hull of the set of measures
$\{\mu_x; x\in C_\xi \quad \mbox{is quasi-regular}\}$; hence
\begin{eqnarray}
\langle \xi, \A_{\mu}\rangle \leq -\frac{\alpha}{2} <0.\label{funny}
\end{eqnarray}

It is well known that
every (finite) positive $\Phi$-invariant Borel measure is supported on $R=R_\xi\cup C_\xi$, see
e.g. \cite{K}, Proposition 4.1.18, page 141. As $R_\xi$ and $C_\xi$ are closed and flow-invariant
we may write $\mu=\mu_1+\mu_2$ where $\mu_1, \mu_2$ are $\Phi$-invariant
and $\mu_1$ is supported on $R_\xi$, while $\mu_2$ is supported on $C_\xi$. It folllows from (\ref{funny}) that
$\langle \xi, \A_{\mu_2}\rangle \leq -\frac{\alpha}{2}\cdot \mu_2(C_\xi).$
Further, we claim that $\langle \xi, \A_{\mu_1}\rangle =0$ for the following
reason. Since $\xi|_{R_\xi}=0$ (as a \v Cech cohomology class),
for any smooth closed 1-form $\omega$ on $M$ representing $\xi$ there exists a smooth function $f$
defined on an open neighborhood of $R_\xi$ such that $\omega=df$ near $R_\xi$.
Then we obtain
\begin{eqnarray}
\qquad\langle \xi, \A_{\mu_1}\rangle = \int_M \iota_V(\omega)d\mu_1 = \int_{R_\xi}
\iota_V(\omega)d\mu_1 = \int_{R_\xi} V(f)d\mu_1=0.\label{funny1}
\end{eqnarray}
The last equality holds since the measure $\mu_1$ is $\Phi$-invariant (see e.g. \cite{Sc}, Theorem on page 277).
Finally, as $\A_\mu = \A_{\mu_1} +\A_{\mu_2}$ we see that $\langle \xi, \A_\mu\rangle \leq -\eta\cdot\mu(C_\xi)$ with $\eta=\alpha/2$
which completes the proof of (II) $\Rightarrow$ (III).

The implication (III)$\Rightarrow$(VI) is obvious.

We are left to show the implication (IV) $\Rightarrow$ (I). Our argument uses the technique
of Schwartzman \cite{Sc}. It is to show that under the conditions (IV) there exists a smooth Lyapunov 1-form
for $(\Phi, R_\xi)$ in the class $\xi$. In a first step we prove that there exists a smooth, closed 1-form $\omega_1$ in the class
$\xi$ so that $\iota_V(\omega_1)<0$ on $C_\xi$.
To this end, denote by $\mathcal D\subset C^0(M)$ the space of functions
$$\mathcal D = \{V(f); f: M\to \R \quad \mbox{is smooth}\}$$
and by ${\mathcal C}^-$ the convex cone in $C^0(M)$ consisting of all functions $f\in C^0(M)$ with
$$f(x)<0 \quad \mbox{for all}\quad x\in C_\xi.$$
As $C_\xi$ is compact, the cone $\mathcal C^-$ is open
in the Banach space $C^0(M)$ of continuous functions on $M$, endowed with the usual supremum norm.
Choose an arbitrary smooth, closed 1-form $\omega$ in the class $\xi$.
Assume that $\mathcal C^- \cap (\iota_V(\omega)+\mathcal D) = \emptyset.$
It then follows from the Hahn - Banach Theorem (cf. \cite{Ru}, page 58) that there exists a continuous linear functional
$\Lambda: C^0(M)\to \R$ so that
$$\Lambda|_{\iota_V(\omega)+\mathcal D} \geq 0\quad \mbox{and}\quad \Lambda|_{\mathcal C^-} <0.$$
Since $\iota_V(\omega)+\mathcal D$ is an affine subspace and $\Lambda$ is bounded on it from below,
we obtain that $\Lambda$ restricted to $\mathcal D$
vanishes. According to the Riesz representation theorem (cf. \cite{La}), there exists a Borel measure $\mu$ on $M$ so that
$$\Lambda(f)=\int_M fd\mu$$
for any $f\in C^0(M)$. By Theorem \cite{Sc}, page 277, the condition $\Lambda|_{\mathcal D}=0$ implies
that $\mu$ is $\Phi$-invariant. On the other hand, $\Lambda|_{\mathcal C^-}<0$ implies that $\mu|_{C_\xi}>0$.

Denote by $\chi: M\to \R$ the characteristic function of $C_\xi$ and let $\nu=\chi\cdot \mu$.
As $C_\xi$ is $\Phi$-invariant $\nu$ is a $\Phi$-invariant Borel measure and (unlike, possibly, $\mu$) is positive.
Note that $\mu -\nu$ is a $\Phi$-invariant Borel measure supported on $R_\xi$ (again using that any $\Phi$-invariant measure
is supported on $R=R_\xi\cup C_\xi$).
Thus, it follows from our assumption $\xi|_{R_\xi}=0$,
by the same argument
which led to (\ref{funny1}),
that
$$\langle \xi, \A_{\mu-\nu}\rangle =0.$$
Since $\A_{\mu-\nu}=\A_\mu -\A_\nu$ we find
$$\langle \xi , \A_\nu\rangle =\langle \xi, \A_\mu\rangle = \int_M \iota_V(\omega)d\mu = \Lambda (f) \geq 0$$
where $f=\iota_V(\omega)$, contradicting condition (IV). This means that the intersection
$\mathcal C^- \cap (\iota_V(\omega)+\mathcal D)$ cannot be empty, i.e. there exists a smooth function $g:M\to \R$
so that the smooth closed 1-form $\omega_1 = \omega+dg$ is in the class $\xi$ and satisfies
$$\iota_V(\omega_1) <0\quad \mbox{on}\quad C_\xi.$$
This completes the first step of the proof.

To finish the argument, we now adjust $\omega_1$ on the complement of
$C_\xi$ so that the resulting form is a Lyapunov 1-form for $(\Phi, R_\xi)$.
As $\iota_V(\omega_1)<0$ on $C_\xi$ and $C_\xi$ is compact, there is some open
neighborhood $W_1$ of $C_\xi$ such that \mbox{$W_1 \cap R_\xi = \emptyset$}
and $\iota_V(\omega_1)<0$ on $W_1$. Since $\xi|_{R_\xi}=0$, there exists an
open neighborhood $W_2$ of $R_\xi$ such that $W_1 \cap W_2 = \emptyset$ and a
smooth function $g:M \to \R$ such that $\omega_{1|W_2} =dg$ and
$dg|_{W_1}=0$.
By Proposition \ref{conley} there exists a smooth Lyapunov function $L: M\to
\R$ for $(\Phi, R)$.
Now consider
\begin{eqnarray}
\omega_2 =\omega_1 - dg +\lambda dL,\label{omegatwo}
\end{eqnarray}
where $\lambda >0$ remains to be chosen. Clearly, the form $\omega_2$ is
smooth and closed and represents the class $\xi$. For any $\lambda>0$
it satisfies $\omega_{2}|_{W_2}=d(\lambda L)$, because $\omega_1 - dg$
vanishes on this set. In particular, $\omega_2$ has property
($\Lambda$2) of a Lyapunov 1-form for the pair $(\Phi,R_\xi)$. Note also
that for all positive $\lambda$ we have $\iota_V(\omega_2)<0$ on $W_1$ (by the
construction of $W_1$) and on $W_2-R_\xi$ because $\omega_1 -dg$ vanishes
there, whereas $V(L)<0$. As the complement of $W_1 \cup W_2$ is compact and
disjoint from $R$,
$$
1 < \lambda_0:= 1 + \sup_{x \notin W_1 \cup W_2}
\frac{|\iota_V(\omega_1-dg)|}{|V(L)|}  <\infty,
$$
and $\iota_V(\omega_2) <0$ on $M-R_\xi$ for all $\lambda \geq \lambda_0$, showing
that for such choices of $\lambda$ the form $\omega_2$ also has property
($\Lambda$1) of a Lyapunov 1-form for $(\Phi,R_\xi)$. This completes the proof of
the implication (IV) $\Rightarrow$ (I) and hence the proof of Theorem
\ref{main}.
\qed

\section*{Appendix: Proof of Proposition \ref{conley}}

Recall from \cite[II 6.2.A]{Conl} the alternative characterization of the chain
recurrent set $R$ as
$$
R = \bigcap\; \{ A \cup A^* \; | \; (A,A^*) \text{ \rm is an
attractor-repeller pair} \}
$$
Here a closed, flow-invariant subset $A\subset M$ is called an attractor if it
admits a neighborhood $U$ such that $A$ is the maximal flow-invariant subset in
the closure of $U \cdot [0,\infty)$. The dual repeller $A^*$ is the set of all
points $x\in M$ whose forward limit set is disjoint from $A$ (cf. \cite[II
5.1]{Conl}). Equivalently, $(A,A^*)$ is an attractor-repeller pair if and only
if both $A$ and $A^*$ are closed flow invariant subsets of $M$ and the forward
(resp. backward) limit set of every point $x\notin A \cup A^*$ is contained in
$A$ (resp. $A^*$) - see \cite[Prop.~1.4.]{RS}.

As $M$ is a closed manifold and hence separable, the number of distinct
attractor-repeller pairs is at most countable (cf. \cite[II 6.4.A]{Conl}). Let
$\{(A_n,A_n^*)\}_{n\geq 1}$ be some enumeration. For each $n\geq 1$, the
construction of Robbin and Salamon (Prop.~1.4. of \cite{RS} and the remark
following it) yields a smooth function \mbox{$f_n:M\to [0,1]$} with
$f_n^{-1}(0)=A_n$, $f_n^{-1}(1)=A^*_n$ and $df_n(V)<0$ on the complement of
$A_n\cup A^*_n$.  Let $c_n$ be positive constants such that in a fixed finite
atlas of charts all partial derivatives of $f_n$ of order $\leq n$ are bounded
pointwise in absolute value by $c_n$. Then
$$
L(x) := \sum_{n=1}^\infty \frac{f_n(x)}{2^nc_n}
$$
is a smooth function having the required properties. In particular, as for any $n\geq 1$
the differential $df_n$ vanishes on $A_n \cup A^\ast_n$, the differential of $L$ vanishes on $R$. \qed

\bibliographystyle{amsalpha}

\begin{thebibliography}{99}


\bibitem{Conl} C. Conley, {\em Isolated invariant sets and the Morse index},
CBMS regional conference series in mathematics, no. 38, AMS, 1976.
\bibitem{Conl1} C. Conley, {\em The gradient structure of a flow: I},
Ergod. Th. \& Dynam. Sys., \textbf{8}(1988), 11--26

\bibitem{D} A. Dold, \textit{Lectures on Algebraic Topology}, Springer - Verlag, 1972.

\bibitem{Fa:1} M. Farber, {\em Zeroes of closed 1-forms, homoclinic orbits
and Lusternik--Schnirelman theory}, "Topological Methods in
Nonlinear Analysis", \textbf{19}(2002), 123 - 152

\bibitem{Fa:2} M. Farber, {\em Lusternik--Schnirelman Theory and Dynamics},
 in: ``Lusternik - Schnirelmann Category and related topics'', Contemporary Mathematics, \textbf{316}(2002).


\bibitem{FKLZ} M. Farber, T. Kappeler, J. Latschev, E. Zehnder,
\textit{Lyapunov 1-forms for flows}, to appear in "Ergodic Theory and Dynamical
Systems", available also as math.DS/0210473


\bibitem{Fr} J. Franks, {\em A Variation on the Poincar\'e-Birkhoff
Theorem}, Contemporary Mathematics, \textbf{81}(1988),  111--117.
\bibitem{Fra:2} J. Franks, {\em Homology and dynamical systems},
CBMS regional conference series in mathematics, no. 49, AMS, 1982.

\bibitem{Fr:1} D. Fried, {\em The geometry of cross sections to
flows}, Topology, vol.~21 (1982), 353--371

\bibitem{Fu:1} F.~B. Fuller, {\em On the surface of section and periodic
trajectories}, Amer. J. Math., \textbf{87} (1965), 473--480

\bibitem{Ja} K. Jacobs, {\em Neuere Methoden in der Ergodentheorie},
Springer, 1960


\bibitem{La} S. Lang, \textit{Real and Functional Analysis}, third edition, Graduate Texts in Mathematics, 142, 1993, Springer - Verlag

\bibitem{K} A. Katok, B. Hasselblatt, \textit{Introduction to the modern theory of dynamical systems},
Cambridge University Press, 1995

\bibitem{RS} J.~Robbin and D.~Salamon, {\em Lyapunov maps, simplicial
complexes and the Stone functor}, Ergod. Th. \& Dynam. Sys. , 12 (1992),
pp. 153--183


\bibitem{Ru} W. Rudin, {\em Functional Analysis}, McGraw-Hill Book Company, 1973

\bibitem{Sc} S. Schwartzman, {\em Asymptotic cycles}, Ann. of Math.,
vol. 66, 1957, 270--284

\bibitem{Sc1} S. Schwartzman, {\em Global cross-sections of compact dynamical systems},
Proc. Nat. Acad. Sci. USA \textbf{48}(1962), 786 - 791


\bibitem{Sh} M. Shub, {\em Global Stability of Dynamical Systems},
Springer Verlag, 1986

\bibitem{Sp} E. Spanier, {\em Algebraic Topology}, Springer Verlag, 1966

\bibitem{Wi} F. W. Wilson, {\em Smoothing derivatives of functions and applications},
Trans. Amer. Math. \textbf{139}(1969), 413 - 428


\end{thebibliography}

\end{document}